\documentclass{sig-alternate}

\usepackage{amsmath,amssymb,stmaryrd,enumerate,bbm,theorem,color}

\usepackage{graphicx}

\newcommand{\assign}{:=}
\newcommand{\tmdummy}{$\mbox{}$}
\newcommand{\tmem}[1]{{\em #1\/}}
\newcommand{\tmmathbf}[1]{\ensuremath{\boldsymbol{#1}}}
\newcommand{\tmop}[1]{\ensuremath{\operatorname{#1}}}
\newcommand{\tmstrong}[1]{\textbf{#1}}

\newenvironment{enumeratealpha}{\begin{enumerate}[a{\textup{)}}] }{\end{enumerate}}
\newenvironment{enumeratenumeric}{\begin{enumerate}[1.] }{\end{enumerate}}
\newenvironment{itemizedot}{\begin{itemize} }{\end{itemize}}

\newenvironment{tmparmod}[3]{\begin{list}{}{\setlength{\topsep}{0pt}\setlength{\leftmargin}{#1}\setlength{\rightmargin}{#2}\setlength{\parindent}{#3}\setlength{\listparindent}{\parindent}\setlength{\itemindent}{\parindent}\setlength{\parsep}{\parskip}} \item[]}{\end{list}}
\definecolor{grey}{rgb}{0.75,0.75,0.75}
\definecolor{orange}{rgb}{1.0,0.5,0.5}
\definecolor{brown}{rgb}{0.5,0.25,0.0}
\definecolor{pink}{rgb}{1.0,0.5,0.5}
\newtheorem{definition}{Definition}
{\theorembodyfont{\rmfamily}}
\newtheorem{Lemma}{Lemma}
\newtheorem{Proposition}{Proposition}
{\theorembodyfont{\rmfamily}\newtheorem{Remark}{Remark}}
\newtheorem{Theorem}{Theorem}

\begin{document}

\conferenceinfo{ISSAC}{'08 Hagenberg, Austria}
\CopyrightYear{2008}

\title{On the Computation of the Topology of a Non-Reduced Implicit
  Space Curve}

\numberofauthors{3}

\author{
\alignauthor
Daouda Niang Diatta\\
\affaddr{University of Limoges, XLIM,}\\
\affaddr{INRIA Sophia-Antipolis, France.}\\
\email{dndiatta@sophia.inria.fr}
\alignauthor
Bernard Mourrain\\
\affaddr{INRIA Sophia-Antipolis, France.}\\
\email{mourrain@sophia.inria.fr}
\alignauthor
Olivier Ruatta\\
\affaddr{University of Limoges, XLIM, France.}\\
\email{olivier.ruatta@unilim.fr}
}

\date{30 July 1999}

\maketitle

\begin{abstract}
  An algorithm is presented for the computation of the topology of a non-reduced space curve defined as the intersection of two implicit algebraic
  surfaces. It computes a Piecewise Linear Structure (PLS) isotopic to the
  original space curve.
  
  The algorithm is designed to provide the exact result for all inputs. It's a
  symbolic-numeric algorithm based on subresultant computation. Simple
  algebraic criteria are given to certify the output of the algorithm.
  
  The algorithm uses only one projection of the non-reduced space curve
  augmented with adjacency information around some ``particular points'' of
  the space curve.
  
  The algorithm is implemented with the Mathemagix Computer Algebra System (CAS) using the SYNAPS library as a backend.
\end{abstract}


\section*{Introduction}
The problem of computing the topological graph of algebraic curves plays an important role in many applications such as plotting  {\cite{MSW}} and sectioning in Computer Aided Geometric Design {\cite{BH}}, {\cite{KCMK}}. A wide literature exists on the computation of the topology of plane
curves ({\cite{GVN}}, {\cite{GLMT}}, {\cite{G2}}, {\cite{OR}}, {\cite{EKW}},
{\cite{El}} and {\cite{H}}). The problem of computing the topology of space
curves has been less investigated. In {\cite{AS}}, Alc\'azar and Sendra give a symbolic-numeric algorithm for {\tmstrong{reduced}} space curves using subresultant and GCD computations of approximated polynomials.  If their approach gives good practical results however it
doesn't give a rigorous proof that a sufficient precision is selected for
all inputs in the computation of GCD of approximated polynomials. In  {\cite
{OR}}, Owen, Rockwood and Alyn give a numerical algorithm for
{\tmstrong{reduced}} space curve using subdivision method. Their algorithm
has a good complexity but the topology around the singularities of the space
curve is not certified. We also mention the work in \ {\cite{GLMT}}, where
two projections of a {\bf reduced} space curve are used, and where the
connection algorithm is valid under genericity conditions.

To our knowledge, the general problem of computing the topology of {\tmstrong{non-reduced}} space curves is not investigated in the algorithmic point of view despite its significance in the problem of computing the topology of a real algebraic surface.
  
We present a certified algorithm that computes the topology of
{\tmstrong{non-reduced}} algebraic space curves.
We compute the topology of a plane projection of the
space curve  and then we lift the computed topology on the space.
The topology of the projected curve is computed using a
classical sweeping algorithm (see {\cite{G2}}, {\cite{GVN}}). For the computation
of the topology of a plane algebraic curve, we present an
{\tmstrong{efficient generic test}} that certifies the output of the algorithm in {\cite{GVN}} .
  
For space curves, we introduce the notion of pseudo-generic position. A space curve
 is said to be in pseudo-generic position with respect to the
$(x, y)$-plane if and only if almost every point of its projection on the $(x,
y)$-plane has only one geometric inverse-image. A simple algebraic criterion is
given to certify the pseudo-genericity of the position of a space curve.
From a theoretical point of view, the use of the notion of curve in
{\tmstrong{pseudo-generic position gives us a rational parametrization of the
space curve}}. The use of this rational parametrization allows us to lift the
topology computed after projection without any supplementary effort.
From a practical point of view, the use of the rational 
parametrization of the space curve makes the lifting faster, avoiding numerical problems.
  
We need to distinguish two kinds of singularities on
the projected curve. A {\tmstrong{certified}} algorithm is given to
do so.\\ Unlike
previous approaches, our algorithm uses {\tmstrong{only one projection}} of the
space curve and works for {\tmstrong{non-reduced space curves.}} We therefore avoid
the cost of the second projection used by previous approaches.\\
 In the next section we describe the fundamental algebraic tools that we use in this paper. In Section
2, we present our contribution to certify the algorithm for computing the
topology of a plane algebraic curve. Our algorithm itself
is introduced in Section 3. We report on our implementation and experiments in section 4.

\section{Subresultants}

Let $P_1, P_2 \in \mathbbm{Q}[X, Y, Z]$ and $\mathcal{C}_{\mathbbm{R}} : =
\{ (x, y, z) \in \mathbbm{R}^3 | \\ P_1 (x, y, z) = P_2 (x, y, z) = 0
\}$ be the intersection of the vanishing sets of $P_1$ and $P_2$. Our
curve analysis needs to compute a plane projection of
$\mathcal{C}_{\mathbbm{R}}$. Subresultant sequences~are a suitable tool to do
it. For the reader's convenience, we recall their definition and relevant
properties. For all the results of this section, we refer to \ {\cite{BR}},
for proofs.

Let $\mathbb{A}$ be a integral domain. Let $P = \sum^p_{i = 0} a_i X^i$ and
$Q = \sum^q_{i = 0} b_i X^i$ be two polynomials with coefficients in
$\mathbb{A}$. We shall always assume $a_p \neq 0$, $b_q \neq 0$ and $p \geqslant
q$.

Let $\mathbbm{P}_r (\mathbb{A})$ be the set of polynomials in
$\mathbb{A}[X]$ of degree not exceeding $r$, always, with the basis (as an
$\mathbb{A}$-module) $1, X, \ldots, X^r$. If $r < 0$, we set $\mathbbm{P}_r
(\mathbb{A}) = 0$ by convention, and we will identify an element $S = s_0 +
\ldots + s_r X^r$ of $\mathbbm{P}_r (\mathbb{A})$ with the row vector $(s_0,
\ldots, s_r)$.

Let $k$ be an integer such that $0 \leqslant k \leqslant q$, and let $\Psi_k$:
\[ \text{$\mathbbm{P}_{q - k - 1} (\mathbb{A}) \times \mathbbm{P}_{p - k - 1}
   (\mathbb{A}) \rightarrow \mathbbm{P}_{p + q - k - 1} (\mathbb{A})$} \]
be the $\mathbb{A}$-linear map defined by $\Psi_k (\tmmathbf{U},
\tmmathbf{V}) = P\tmmathbf{U}+ Q\tmmathbf{V}$, with $M_k (P, Q)$ the $(p + q -
k) \times (p + q - k)$ matrix of $\Psi_k$. As we write vectors as row vectors,
we have
\[ M_k (P, Q) = \left(\begin{array}{ccccc}
     a_0 & \ldots & a_p &  & \\
     & \ddots &  & \ddots & \\
     &  & a_0 & \ldots & a_p\\
     b_0 & \ldots & b_q &  & \\
     & \ddots &  & \ddots & \\
     &  & b_0 & \ldots & b_q
   \end{array}\right) \]
That is $M_0 (P, Q)$ is the classical Sylvester matrix associated to $P, Q$.
To be coherent with the degree of polynomials, we will attach index $i - 1$ to
the $i^{\tmop{th}}$ column of $M_k (P, Q)$, so the indices of the columns go
from 0 to $p + q - k - 1$.

\begin{definition}
  For $j \leqslant p + q - k - 1$ and $0 \leqslant k \leqslant q$, let $\tmop{sr}_{k, j}$ be the determinant of the submatrix of $M_k (P, Q)$
  formed by the last $p + q - 2 k - 1$ columns, the column of index j and all
  the $(p + q - 2 k)$ rows. The polynomial $\tmop{Sr}_k (P, Q) = \tmop{sr}_{k,
  0} + \ldots + \tmop{sr}_{k, k} X^k$ is the $k^{\tmop{th}}$ sub-GCD of P and
  Q, and its leading term $\tmop{sr}_{k, k}$ $($sometimes noted $\tmop{sr}_k
  )$ is the $k^{\tmop{th}}$ subresultant of P and Q. So, it follows that
  $\tmop{Sr}_0 (P, Q) = \tmop{sr}_0$ is the usual resultant of P and Q.
\end{definition}

\begin{Remark}
  \label{remsubres}{\tmdummy}
  
  \begin{enumeratenumeric}
    \item For $k < j \leqslant p + q - k - 1$, we have $\tmop{sr}_{k, j} = 0$,
    because it is the determinant of a matrix with two equal columns.
    
    \item If $q < p$, we have $\tmop{Sr}_q = (b_q)^{p - q - 1} Q$ and
    $\tmop{sr}_q = (b_q)^{p - q}$.
  \end{enumeratenumeric}
\end{Remark}

The following proposition will justify the name of sub-GCD given to the
polynomial $\tmop{Sr}_k$.

\begin{Proposition}
  \label{propsubres}Let d be the degree of the GCD of P and Q $($d is defined
  because $\mathbb{A}$ is an integral domain, so we may compute the GCD over
  the quotient field of $\mathbb{A} )$. Let k be an integer such that $k
  \leqslant d$.
  \begin{enumeratenumeric}
    \item The following assertions are equivalent:
    \begin{enumeratealpha}
      \item $k < d$;
      
      \item $\tmop{Sr}_k = 0$;
      
      \item $\tmop{sr}_k = 0$.
    \end{enumeratealpha}
    \item $\tmop{sr}_d \neq 0$ and $\tmop{Sr}_d$ is the GCD of
    P and Q.
  \end{enumeratenumeric}
\end{Proposition}

\begin{Theorem}
  \label{thmsubres}{\tmstrong{Fundamental property of subresultants}} \\
  The first polynomial $\tmop{Sr}_k$ associated to P and Q with $\tmop{sr}_k
  \neq 0$ is the greatest common divisor of P and Q.
\end{Theorem}

We will often call $(\tmop{Sr}_i)_i$ the subresultant sequence associated to
$P$ and $Q$ and $(\tmop{sr}_{i, j})_{i, j}$ the sequence of their subresultants
coefficients.
We will denote by $\tmop{lcoef}_X (f)$ the leading coefficient of the polynomial $f$ with respect to the variable $X$.

\begin{Theorem}
  \label{thmsubres}{\tmstrong{Specialization \ property \ of
  subresultants}} \\
  Let $P_1, P_2 \in \mathbb{A}[Y, Z]$ and $\left( \tmop{Sr}_i (Y, Z) \left)_i
  \right. \right.$ be their subresultant sequence with respect to $Z$. Then
  for any $\alpha \in \mathbb{A}$ with: \\
 deg$_Z (P (Y, Z))$ = deg$_Z (P (\alpha,Z))$;\\
 deg$_Z (Q (Y, Z))$ = deg$_Z (Q (\alpha, Z))$,\\
  $\left( \tmop{Sr}_i(\alpha, Z) \left)_i \right. \right.$ is the subresultant sequence of the
  polynomials $P (\alpha, Z)$ and $Q (\alpha, Z)$.
\end{Theorem}

\section{Topology of a plane algebraic curve}

Let $f \in \mathbbm{Q}[X, Y]$ be a square free polynomial and \begin{equation} \mathcal{C}(f)
\assign \{ (\alpha, \beta) \in \mathbbm{R}^2, f (\alpha, \beta) = 0 \} \end {equation} be the
real algebraic curve associated to $f$. We want to compute the topology of
$\mathcal{C}(f)$.

For curves in generic position, computing its critical fibers and one regular fiber between two
critical ones is sufficient to obtain the topology using a sweeping
algorithm (see {\cite{GVN}}). But for a good computational behaviour, it is
essential to certify the genericity of the position of the curve.
  
We propose an effective test allowing to certify the computation and connection, in a deterministic way. This is an important tool in
order to address the case of space curves.

Now, let us introduce the definitions of generic position, critical, singular
and regular points.

\begin{definition}
  \label{defxcrit}Let $f \in \mathbbm{Q}[X, Y]$ be a square free polynomial
  and $\mathcal{C}(f) = \{ (\alpha, \beta) \in \mathbbm{R}^2 : f (\alpha,
  \beta) = 0 \}$ be the curve defined by f. A point $(\alpha, \beta) \in
  \mathcal{C}(f)$ is called:
  \begin{itemizedot}
   \item a x-critical point if $\partial_Y f (\alpha,
    \beta) = 0$,
    
   \item a singular point if $\partial_X f (\alpha,
    \beta) = \partial_Y f (\alpha, \beta) = 0$,
    
   \item a regular point if $\partial_X f (\alpha,\beta) \neq 0$ or $\partial_Y f (\alpha, \beta) \neq 0$.
  \end{itemizedot}
\end{definition}

With these definitions we can describe the generic conditions required for
plane curves.

\begin{definition}
  \label{defgenpos}Let $f \in \mathbbm{Q}[X, Y]$ be a square free polynomial
  and $\mathcal{C}(f) = \{ (\alpha, \beta) \in \mathbbm{R}^2 : f (\alpha,
  \beta) = 0 \}$ be the curve defined by f. Let $\mathcal{N}_x (\alpha) :=\# \{
  \beta \in \mathbbm{R}$, such that $(\alpha, \beta)$ is a x-critical point of
  $\mathcal{C}(f) \left\} \right.$. $\mathcal{C}(f)$ is in generic position
  for the x-direction, if:
  \begin{enumeratenumeric}
    \item $\forall \alpha \in \mathbbm{C}, \mathcal{N}_x (\alpha) \leqslant
    1,$
    
    \item There is no asymptotic direction of $\mathcal{C}(f)$ parallel to the
    y-axis. 
  \end{enumeratenumeric}
\end{definition}

This notion of genericity also appears in \ {\cite{EKW}} and in a slightly
more restrictive form in \ {\cite{El}}. Previous approaches succeed if
genericity conditions are satisfied, but they do {\tmstrong{not guarantee to reject the curve}} if
they are not; i.e, it does {\tmstrong{not decide genericity}}. So
for some input curves the computed topology might not be exact.

A change of coordinates such that lcoef$_Y (f) \in \mathbbm{Q}^{\ast}$ is
sufficient to place $\mathcal{C}(f)$ in a position such that any asymptotic
direction is not parallel to the $y$-axis. It remains to find an efficient way to verify \ the first condition. This
follows from the next propositions. We refer to {\cite{GVN}}, for proofs.

\begin{Proposition}
  Let $f \in \mathbbm{Q} \left[ X, Y] \right.$ be a square free polynomial
  with $\tmop{lcoef}_Y (f) \in \mathbbm{Q}^{\ast}$, $\tmop{Res}_Y (f,
  \partial_Y f)$ be the resultant with respect to Y of the polynomials $f$,
  $\partial_Y f$ and $\left. \right\{ \alpha_1, \ldots, \alpha_l \}$ be the set
  of the roots of $\tmop{Res}_Y (f, \partial_Y f)$ in $\mathbbm{C}$.
  
  Then $\mathcal{C}(f)$ is in generic position if and only if\\
 $\forall i \in\{1, \ldots, l\},  \gcd \left( f (\alpha_i, Y), \partial_Y f (\alpha_i, Y)
  \right)$ has at most one root.
\end{Proposition} \label{prop:genposcarac}

Let $f \in \mathbbm{Q} \left[ X, Y] \right.$ be a square free polynomial with\\
$\tmop{lcoef}_Y (f) \in \mathbbm{Q}^{\ast}$ and $d \assign \deg_Y (f)$. We
denote by $\tmop{Sr}_i (X, Y)$ the $i^{\text{th}}$ subresultant polynomial of $f$
and $\partial_Y f$ and $\tmop{sr}_{i, j} (X)$ the coefficient of $Y^j$ in
$\tmop{Sr}_i (X, Y)$. We define \ inductively the following polynomials:
\[ \Phi_0 (X) = \frac{\tmop{sr}_{0, 0} (X)}{\gcd (\tmop{sr}_{0, 0} (X),
   \tmop{sr}'_{0, 0} (X))} ; \]
 $\forall i \in \{ 1, \ldots, d - 1 \},$
$\Phi_i (X) = \gcd(\Phi_{i - 1} (X), \tmop{sr}_{i, i_{}} (X))$ and $\Gamma_i (X) = \frac{\Phi_{i -1} (X)}{\Phi_i (X)}$.

\begin{Proposition} \label{prop:princ}
  {\tmdummy}
  
  \begin{enumeratenumeric}
    \item $\Phi_0 (X) = \underset{i = 1}{\overset{d - 1}{\prod}} \Gamma_i (X)$
    and $\forall i$, $j \in \{1, \ldots, d - 1\}, i \neq j \Longrightarrow
    \gcd (\Gamma_i (X), \Gamma_j (X)) = 1 ;$
    
    \item Let $k \in \{1, \ldots, d - 1\}$, $\alpha \in \mathbbm{C}$.
    $\Gamma_k (\alpha) = 0$ if and only if $\gcd (f (\alpha, Y), \partial_Y f
    (\alpha, Y)) = \tmop{Sr}_k (\alpha, Y)$;
    
    \item $ \{ (\alpha, \beta) \in \mathbbm{R}^2 : f (\alpha, \beta) = \partial_Y f (\alpha, \beta) = 0 \} =
\bigcup^{d - 1}_{k = 1} \{(\alpha, \beta) \in \mathbbm{R}^2 : \Gamma_k (\alpha) = \tmop{Sr}_k
  (\alpha, \beta) = 0  \} .$
  \end{enumeratenumeric}
\end{Proposition}

In the following theorem, we give an effective and efficient algebraic test to
certify the genericity of the position of a curve with respect to a given
direction.

\begin{Theorem} \label{thmtestgenpos}
  Let $f \in \mathbbm{Q} \left[ X, Y] \right.$ be a square free polynomial
  such that deg$_Y (f_{}) = d$, lcoef$_Y (f) \in \mathbbm{Q}^{\ast}$. Then
  $\mathcal{C}(f)$ is in generic position for the projection on the $x$ axis
  if and only if 
   $\forall k \in \{ 1,\ldots, d - 1\}$,
   $\forall i \in \{ 0,\ldots, k - 1 \},$ \\
   $k (k - i) \tmop{sr}_{k, i} (X) \tmop{sr}^{_{}}_{k, k} (X) - (i +
     1) \tmop{sr}_{k, k - 1} (X) \tmop{sr}_{k, i + 1} (X) =$\\
     $0 \tmop{mod}\Gamma_k (X) .$
\end{Theorem}

\begin{proof}
  Assume that $\mathcal{C}(f)$ is in generic position and let $\alpha
  \in \mathbbm{C}$ be a root of $\Gamma_k (X)$. According to Proposition \ref{prop:princ} (2.) \\ $\gcd ( f
  (\alpha, Y), \partial_Y f (\alpha, Y)) = \tmop{Sr}_k (\alpha, Y) =
  \sum_{j = 0}^k \tmop{sr}_{k, j} (\alpha) Y^j
  .$\\ According to Proposition \ref{prop:genposcarac},
  $\tmop{Sr}_k (\alpha, Y)$ has an only root \\$\beta (\alpha) = -
  \frac{\tmop{sr}_{k, k} (\alpha)}{k \tmop{sr}_{k, k - 1} (\alpha)}$,
  so  $\tmop{Sr}_k (\alpha, Y) = \tmop{sr}_{k, k} (\alpha)
  (Y - \beta)^k$.\\Binomial Newton formula gives\\$\tmop{Sr}_k (\alpha,
  Y) = \tmop{sr}_{k, k} (\alpha) (Y - \beta)^k = \tmop{sr}_{k, k}
  (\alpha) \sum_{i = 0}^k \binom{k}{i} (- \beta)^{k - i} Y^i$.
  So by identification $\forall k \in \{1, \ldots, d - 1\}, \forall i \in \{0, \ldots, k -
  1\}$ and $\forall \alpha \in \mathbbm{C}$ such that $\Gamma_k (\alpha) = 0$,
  \[ k (k - i) \tmop{sr}_{k, i} (\alpha) \tmop{sr}_{k, k} (\alpha) - (i + 1)
     \tmop{sr}_{k, k - 1} (\alpha) \tmop{sr}_{k, i + 1} (\alpha) = 0. \]
  It is to say that $\forall k \in \{1, \ldots, d - 1\}, \forall i \in \{0,
  \ldots, k - 1\},$\\
  $k (k - i) \tmop{sr}_{k, i} (X) \tmop{sr}^{_{}}_{k, k} (X) - (i +
     1) \tmop{sr}_{k, k - 1} (X) \tmop{sr}_{k, i + 1} (X) = 0 \tmop{mod}
     \Gamma_k (X) .$\\
  Conversely, let $\alpha$ be a root of $\Gamma_k (X)$ such that \[ k (k - i) \tmop{sr}_{k, i} (\alpha) \tmop{sr}_{k, k} (\alpha) - (i + 1)
     \tmop{sr}_{k, k - 1} (\alpha) \tmop{sr}_{k, i + 1} (\alpha) = 0. \] With the same argument used in the first part of this proof we
  obtain 
\begin{align}
\begin{array}{lcl}
\gcd \left( f (\alpha, Y), \partial_Y f (\alpha, Y) \right) 
& = & \tmop{Sr}_k (\alpha, Y) \\
& = & \sum^k_{j = 0} \tmop{sr}_{k, j} (\alpha) Y^j \\
& = & \tmop{sr}_{k, k} (\alpha) (Y - \beta)^k 
\end{array}
\end{align}
\noindent with 
\begin{align} 
\label {betaa}
\beta = - \frac{\tmop{sr}_{k, k - 1 (\alpha)}}{k \tmop{sr}_{k, k (\alpha)}}.
\end{align} 
Then we conclude that \ $\left. \gcd ( f (\alpha_{}, Y), \partial_Y f
(\alpha_{}, Y) \right)$ has only one distinct root and, according to
Proposition \ref{prop:genposcarac}, $\mathcal{C}(f)$ is in generic
position.
\end{proof}

\begin{Remark}
  \label{remgenpos}Theorem \ref{thmtestgenpos} shows that it is possible to check
  with certainty if a plane algebraic curve is in generic position or not. If not, we can put it in generic position by a basis
  change.\\
  In fact, it is well known that there is only a finite number of bad changes
  of coordinates of the form $X : = X + \lambda Y$, $Y : = Y$,
  such that if $\mathcal{C}(f)$ is not in generic position then the
  transformed curve remains in a non-generic position. This number of bad
  cases is bounded by $\binom{c}{2}$, where $c$ is the number of distinct $x$-critical points of $\mathcal{C}(f)$ \ {\cite{GVN}}.
\end{Remark}

\section{Topology of implicit three dimensional algebraic curves}

\subsection{Description of the problem }

Let $P_1, P_2 \in \mathbbm{Q}[X, Y, Z]$ and \begin{equation}\mathcal{C}_{\mathbbm{R}}:=
\{ (x, y, z) \in \mathbbm{R}^3:\\ P_1 (x, y, z) = P_2 (x, y, z) = 0
\} \end{equation} be the intersection of the surfaces defined by $P_1 = 0$ and $P_2 =
0$. We assume that $\gcd (P_1, P_2) = 1$ so that $\mathcal{C}_{\mathbbm{R}}$
is a space curve. Our goal is to analyze the geometry of
$\mathcal{C}_{\mathbbm{R}}$ in the following sense: We want to compute a
piecewise linear graph of $\mathbbm{R}^3$ isotopic to the original space
curve.\\
Our method allows to use a new sweeping algorithm using only one projection of
the space curve.\\
To make the lifting possible using only one projection, \ a new definition
of generic position for space curves and an algebraic characterization of it
are given. We will also need to distinguish the "apparent singularities" and
the "real singularities". A {\tmstrong{certified}} algorithm is given to
{\tmstrong{distinguish}} these two kinds of {\tmstrong{singularities}}.\\
For the lifting phase, using the new notion of curve in
{\tmstrong{pseudo-generic position}}, we give an algorithm that computes a
rational parametrization of the space curve. The use of this rational
parametrization allows us to lift the topology of the projected curve
without any supplementary computation.

\subsection{Genericity conditions for space curves}

Let $\Pi_z : (x, y, z) \in \mathbbm{R}^3 \mapsto (x, y) \in \mathbbm{R}^2$. We
still denote $\Pi_z = \Pi_z |_{\mathcal{C}_{\mathbbm{R}}}$. Let $\mathcal{D}=
\Pi_z (\mathcal{C}_{\mathbbm{R}}) \subset \mathbbm{R}^2$ be the curve obtained by
projection of $\mathcal{C}_{\mathbbm{R}}$.
  
We assume that $\deg_Z \left( P_1 \right) = \deg (P_1)$ and $\deg_Z (P_2)
= \deg (P_2)$ (by a basis change, these conditions are always satisfied). Let $h (X, Y)$ be {\tmstrong{the squarefree}} part of $\tmop{Res}_Z (P_1, P_2)
\in \mathbbm{Q}[X, Y]$. With the above notation and assumptions we have the
following "geometric" equality, $\Pi_z (\mathcal{C}_{\mathbbm{R}})
=\mathcal{C}(h)$.

\begin{definition}
  \label{defpseudogenpos}{\tmstrong{[Pseudo-generic position]}}
  
  Let $\mathcal{C}_{\mathbbm{C}} : = \left\{ (x, y, z) \in \mathbbm{C}^3 | P_1
  (x, y, z) = P_2 (x, y, z) = 0 \right\}$.\\ The curve
  $\mathcal{C}_{\mathbbm{R}}$ is in pseudo-generic position with respect to
  the $(x, y)$-plane if and only if almost every point of $\Pi_z
  (\mathcal{C}_{\mathbbm{C}})$ has only one geometric inverse-image, i.e.
  generically, if $(\alpha, \beta) \in \Pi_z (\mathcal{C}_{\mathbbm{C}})$,
  then $\Pi^{- 1}_z (\alpha, \beta)$ consists in one point possibly multiple.
\end{definition}
Let $m$ be the minimum of $\deg_Z \left( P_1 \right)$ and $\deg_Z (P_2)$.\\
The following theorems give us an effective way to test if a curve is in
pseudo-generic position or not.

\begin{Theorem}\label{thmpseudogenpos}Let $( \tmop{Sr}_j (X, Y, Z))_{_{j \in \{ 0,\ldots,
  m \left\} \right.}}$ be the subresultant sequence and $\left(
  \tmop{sr}_j (X, Y) \right)_{j \in \{ 0,\ldots m \}}$ be the
  principal subresultant coefficient sequence. Let $( \Delta_i (X, Y))_{i \in
  \left\{ 1,\ldots, m \} \right.}$ be the sequence of $\mathbbm{Q}[X,
  Y]$ defined by the following relations
  \begin{itemizedot}
    \item $\Delta_0 (X, Y) = 1; \Theta_0 (X, Y) = h (X, Y);$
    
    \item For $i \in \{1,...,m \},\\ \Theta_i (X, Y) = \gcd (
    \Theta_{i - 1} (X, Y), \tmop{sr}_i (X, Y)),\\ \Delta_i (X, Y) =
    \frac{\Theta_{i - 1} (X, Y)}{\Theta_i (X, Y)} .$
  \end{itemizedot}
  For $i \in \{ 1,\ldots, m \}$, let $\mathcal{C}(\Delta_i) :=
  \left\{ (x, y) \in \mathbbm{R}^2 | \text{$\Delta_i (x, y) = 0 \left\}
  \right.$} \right.$ and $\mathcal{C}(h) := \{ (x, y) \in \mathbbm{R}^2 |h (x,
  y) = 0 \}$ then
  \begin{enumeratenumeric}
    \item $h (X, Y) = \underset{i = 1}{\overset{m}{\prod}} \Delta_i (X, Y),$
    
    \item $\mathcal{C}(h) = \underset{i = 1}{\overset{m}{\bigcup}}
    \mathcal{C}(\Delta_i),$
    
    \item $\mathcal{C}_{\mathbbm{R}}$ is in pseudo-generic position with
    respect to the $(x, y)$-plane if and only if $\forall i \in \{1, \ldots,
    m\}, \forall (x, y) \in \mathbbm{C}^2$ such that $\tmop{sr}_i (x, y) \neq
    0$ and $\Delta_i (x, y) = 0$, we have \\$\tmop{Sr}_i (x, y, Z) =
    \tmop{sr}_{i, i} (x, y) \left( Z + \frac{\tmop{sr}_{i, i - 1} (x, y)}{i
    \tmop{sr}_{i, i} (x, y)} \right)^i .$
  \end{enumeratenumeric}
\end{Theorem}

\begin{proof}
\begin{enumeratenumeric}
  \item By definition, \ $\forall i \in \{1, \ldots, m\},\\ \Delta_i (X, Y) =
  \frac{\Theta_{i - 1} (X, Y)}{\Theta_i (X, Y)}$. So by a trivial induction
  \[ \text{$\underset{i = 1}{\overset{m}{\prod}} \text{$\Delta_i (X, Y )$} =
     \frac{\Theta_0 (X, Y)}{\Theta_m (X, Y)}$.}  \]
  $\deg_Z (P_1) = \deg (P_1)$ and $\deg_Z (P_2)
  = \deg (P_2)$ imply\\$\tmop{sr}_m (X, Y)
  \in \mathbbm{Q}^{\ast}$ (see $\tmop{Remark} \ref{remsubres}$).\\So  $\Theta_m (X, Y) = \gcd ( \Theta_{m - 1} (X, Y),
  \tmop{sr}_m (X, Y) \left) = 1,\right.$ then\\$\underset{i = 1}{\overset{m}{\prod}}
  \text{$\Delta_i (X, Y )$} = \Theta_0 (X, Y) = h (X, Y)$.

  \item Knowing that $h (X, Y) = \underset{i = 1}{\overset{m}{\prod}}
  \text{$\Delta_i (X, Y )$}$, so it is clear that $\mathcal{C}(h) = \underset{i
  = 1}{\overset{m}{\bigcup}} \mathcal{C}({\Delta_i})$.

  \item Assume that $\mathcal{C}_{\mathbbm{R}}$ is in pseudo-generic position with
  respect to the $( x, y )$-plane. Let $i \in \{1, \ldots, m\}$
  and $(\alpha, \beta)$ $\in \mathbbm{C}^2$ such that sr$_i (\alpha,
  \beta) \neq 0$ and $\Delta_i (\alpha, \beta ) = 0$. Then
  $\Delta_i (X, Y$)$= \frac{\Theta_{i - 1} (X, Y)}{\Theta_i (X, Y)}
   \Longrightarrow$ $\Theta_{i - 1} (\alpha, \beta) = 0$. Knowing that
 $\Theta_{i - 1} (X,Y) = \gcd ( \Theta_{i - 2} (X, Y), \tmop{sr}_{i - 1} (X, Y))$, so it exists
  $d_1, d_2 \in \mathbbm{Q}[X, Y]$ such that \\$\Theta_{i - 2} (X, Y) = d_1 (X,Y) \Theta_{i - 1} (X, Y)$ and\\$\tmop{sr}_{i - 1} (X, Y) = d_2 (X, Y)
  \text{$\Theta_{i - 1} (X, Y)$}$. In this way,\\$\Theta_{i - 1} (\alpha,
  \beta) = 0 \Longrightarrow$ $\Theta_{i - 2} (\alpha, \beta) = 0$ and
  $\tmop{sr}_{i - 1} (\alpha, \beta) = 0$. By using the same arguments, $\Theta_{i - 2}
  (\alpha, \beta) = 0 \Longrightarrow$ $\Theta_{i - 3}
  (\alpha, \beta) = 0$ and $\tmop{sr}_{i - 2} (\alpha, \beta) = 0$. By
  repeating the same argument, we show $\tmop{sr}_{i - 1} (\alpha, \beta)
  = \ldots = \tmop{sr}_0 (\alpha, \beta) = 0.$ Because $\tmop{sr}_i
  (\alpha, \beta) \neq 0$, then the fundamental theorem of subresultant gives\\ $\gcd( (P_1 (\alpha, \beta, Z), P_2 (\alpha, \beta,
  Z)) = \tmop{Sr}_i (\alpha, \beta, Z) =\\ \sum_{j = 0}^i \tmop{sr}_{i, i - j}
  (\alpha, \beta) Z^{i - j} .$ Knowing that \
  $\mathcal{C}_{\mathbbm{R}}$ is in pseudo-generic position with respect to
  the $(x, y)$-plane and\\$\Delta_i (\alpha, \beta) = 0$ then the polynomial
  $\tmop{Sr}_i (\alpha, \beta, Z)$ has only one distinct root which can be
  written $- \frac{\tmop{sr}_{i, i - 1} (\alpha, \beta)}{i \tmop{sr}_{i, i -
  1} (\alpha, \beta)}$
 depending on the relation between coefficients and
  roots of a polynomial. So $\tmop{Sr}_i (\alpha, \beta, Z) = \underset{j =
  0}{\overset{m}{\sum}} \tmop{sr}_{, i, i - j} (\alpha, \beta) Z^{i - j} =\\
  \tmop{sr}_{i, i} (\alpha, \beta) \left( Z + \frac{\tmop{sr}_{i, i - 1}
  (\alpha, \beta)}{i\tmop{sr}_{i, i - 1} (\alpha, \beta)} \right)^i .$\\
  Conversely, assume that $\forall i \in \{1, \ldots, m\}$, $\forall (x, y)
  \in \mathbbm{C}^2$ such that $\tmop{sr}_i (x, y) \neq 0$ and $\Delta_i (x,
  y) = 0,$ we have\\
 $\tmop{Sr}_i (x, y, Z) = \underset{j =0}{\overset{m}{\sum}} \tmop{sr}_{, i, i - j} (x, y) Z^{i - j} =\\
  \tmop{sr}_{i, i} (x, y) \left( Z + \frac{\tmop{sr}_{i, i - 1} (x, y)}{i
  \tmop{sr}_{i, i - 1} (x, y)} \right)^i .$ Let $\mathcal{O}$ be an irreductible component
  of $\Pi_z (\mathcal{C}_{\mathbbm{C}})$. Then there exists $i \in \{1,
  \ldots, m\}$ such that $\mathcal{O} \subset \mathcal{C}(\Delta_i)$. Let
  $(\alpha, \beta) \text{}$ be a point of $\mathcal{O}$, such that $\Delta_i (\alpha, \beta) = 0$ and
  $\tmop{sr}_i (\alpha, \beta) \neq 0$. Now if we define $\gamma \assign -
  \frac{\tmop{sr}_{i, i - 1} (\alpha, \beta)}{i \tmop{sr}_{i, i} (\alpha,
  \beta)}$, we obtain that $\tmop{Sr}_i (\alpha, \beta, \gamma) = 0$, then \ $(\alpha,
  \beta, \gamma)$ is the only point of $\mathcal{C}_{\mathbbm{C}}$ with
  $(\alpha, \beta)$ as projection. So $\mathcal{C}_{\mathbbm{R}}$ is in
  pseudo-generic position with respect to the $\left( x, y \right)$-plane.
\end{enumeratenumeric}
\end{proof}
The following proposition is a corollary of the third result of the previous theorem. If $\mathcal{C}_{\mathbbm{R}}$ is in pseudo-generic position with
  respect to the $(x, y)$-plane, it gives a rational parametrization for
  the regular points of $\mathcal{C}_{\mathbbm{R}}$.

\begin{Proposition} \label{proplifting}
Assume that $\mathcal{C}_{\mathbbm{R}}$ is in pseudo-generic position with
  respect to the $( x, y )$-plane and let \ $(\alpha,
  \beta, \gamma)\in \mathcal{C}_{\mathbbm{R}}$ such that $\tmop{sr}_i (\alpha, \beta) \neq
    0$ and $\Delta_i (\alpha, \beta) = 0$. Then, \begin{equation} \gamma \assign -
  \frac{\tmop{sr}_{i, i - 1} (\alpha, \beta)}{i \tmop{sr}_{i, i} (\alpha,
  \beta)}. \end{equation}
\end{Proposition}
\begin{Remark}
By construction, the parametrization given in Proposition \ref{proplifting} is valid when $\tmop{sr}_{i, i} (\alpha,
  \beta) \neq 0$. If $\tmop{sr}_{i, i} (\alpha,\beta) = 0$ then either $\Delta_j (\alpha, \beta) = 0$ for some $j>i$ or $(\alpha, \beta)$ is a $x$-critical point of $\mathcal{C}(\Delta_i)$ (see section 3.3).
\end {Remark}
The following theorem gives an algebraic test to
certify the pseudo-genericity of the position of a space curve with respect to
a given plane.

\begin{Theorem}
  \label{thm3dgenposcheck}Let $\left( \tmop{Sr}_j (X, Y, Z) \right)_{j \in
  \{0, \ldots, m\}}$ be the subresultants sequence associated to $P_1 (X, Y,
  Z)$ and $P_2 (X, Y, Z)$ and\\$(\Delta_i (X, Y))_{i \in \{1, \ldots, m\}}$
  be the sequence of $\mathbbm{Q}[X, Y]$ previously defined. The curve
  $\mathcal{C}_{\mathbbm{R}}$ is in pseudo-generic position with respect to
  the $(x, y)$-plane if and only if
  \[ \forall i \in \{1, \ldots, m - 1\}, \forall j \in \{0, \ldots, i - 1\},
  \]
   $i (i - j) \tmop{sr}_{i, j} (X, Y) \tmop{sr}_{i, i} (X, Y) - (j + 1)
     \tmop{sr}_{i, i - 1} (X, Y)\\ \tmop{sr}_{i, j + 1} (X, Y) = 0 \tmop{mod}
     \Delta_i (X, Y).$
\end{Theorem}

\begin{proof}
  Assume $\mathcal{C}_{\mathbbm{R}}$ be in pseudo-generic position. Let $i \in
  \{1, \ldots, m - 1\}$, $j_{} \in \{0, \ldots, i - 1\}$, $(\alpha, \beta) \in
  \mathbbm{R}^2$ such that \ \ \ $\Delta_i (\alpha, \beta) = 0$.\\
  If $\tmop{sr}_{i, i} (\alpha, \beta) = 0,$ then $\tmop{sr}_{i, i - 1}
  (\alpha, \beta) = 0$, consequently
 $i (j + 1) \tmop{sr}_{i, j + 1} (\alpha,
  \beta) \tmop{sr}_{i, i} (\alpha, \beta) - (i - j) \tmop{sr}_{i, i - 1}
  (\alpha, \beta) \tmop{sr}_{i, j} (\alpha, \beta) = 0.$\\
  If $\tmop{sr}_{i, i} (\alpha, \beta) \neq 0,$ then according to Theorem \ref{thmpseudogenpos} (3.)\\
 $\tmop{Sr}_i (\alpha, \beta, Z) =
  \underset{j = 0}{\overset{i}{\sum}} \tmop{sr}_{, i, i - j} (\alpha, \beta)
  Z^{i - j} = \\\tmop{sr}_{i, i} (\alpha, \beta) \left( Z + \frac{\tmop{sr}_{i,
  i - 1} (\alpha, \beta)}{i \tmop{sr}_{i, i - 1} (\alpha, \beta)} \right)^i .$
  Let $\gamma \assign - \frac{\tmop{sr}_{i, i - 1} (\alpha, \beta)}{i
  \tmop{sr}_{i, i - 1} (\alpha, \beta)}$, then $\tmop{Sr}_i (\alpha,
  \beta, Z) = \underset{j = 0}{\overset{i}{\sum}} \tmop{sr}_{, i, i - j}
  (\alpha, \beta) Z^{i - j} = \tmop{sr}_{i, i} (\alpha, \beta) \left( Z -
  \gamma \right)^i .$ Using the binomial Newton
  formula we obtain
 $\tmop{Sr}_i (\alpha, \beta,
  Z) = \\\underset{j = 0}{\overset{i}{\sum}} \tmop{sr}_{, i, i - j} (\alpha,
  \beta) Z^{i - j} = \tmop{sr}_{i, i} (\alpha, \beta) \underset{j =
  0}{\overset{i}{\sum}} \binom{i}{j} (- \gamma)^{i - j} Z^j .$ So by
  identification, it comes that
  \[ \forall i \in \{1, \ldots, m - 1\}, \forall j \in \{0, \ldots, i - 1\},
  \]
   $i (i - j) \tmop{sr}_{i, j} (\alpha, \beta) \tmop{sr}_{i, i}
     (\alpha, \beta) - (i + j) \tmop{sr}_{i, i - 1} (\alpha, \beta)
     \tmop{sr}_{i, j + 1} (\alpha, \beta) = 0$, $\forall (\alpha,
     \beta),  \Delta_i (\alpha, \beta) = 0$.
  The reciprocal uses the same arguments.
\end{proof}

\begin{Remark}
  Theorem \ref{thm3dgenposcheck} shows that it is possible to check with
  certainty if a space algebraic curve is in pseudo-generic position or not.
  If it is not, we can put it in pseudo-generic position by a change of
  coordinates. In fact, there is only a finite number of bad changes of
  coordinates of the form
  
  $X : = X + \lambda Z$; $Y : = Y + \mu Z$; $Z \assign Z$,\\
  with $\lambda, \mu \in \mathbbm{Q}^{\ast}$ such that if
  $\mathcal{C}_{\mathbbm{R}}$ is not in pseudo-generic position then the
  transformed curve remains in a non-pseudo-generic position {\cite{AS}}.
\end{Remark}
Let us introduce the definitions of generic position, critical, singular,
regular points, apparent singularity and real singularity for a space
algebraic curve.

\begin{definition}
  Let $M(X, Y, Z)$ be the 2$\times$3 Jacobian matrix with rows $(\partial_X P_1,
  \partial_Y P_1, \partial_Z P_1)$ and $(\partial_X P_2, \partial_Y P_2,
  \partial_Z P_2)$.
  \begin{itemizedot}
    \item A point p$\in \mathcal{C}_{\mathbbm{R}}$ is regular (or smooth) if
    the rank of $M(p)$ is 2.
    
    \item A point p$\in \mathcal{C}_{\mathbbm{R}}$ which is not regular is
    called singular.
    
    \item A point $p = (\alpha, \beta, \gamma) \in \mathcal{C}_{\mathbbm{R}}$
    is x-critical (or critical for the projection on the x-axis) if the curve
    $\mathcal{C}_{\mathbbm{R}^{}}$ is tangent at this point to a plane
    parallel to the $($y,z$)$-plane. The corresponding $\alpha$ is called a
    x-critical value.
  \end{itemizedot}
\end{definition}

\begin{figure}
  \centering
  \includegraphics[width=\linewidth]{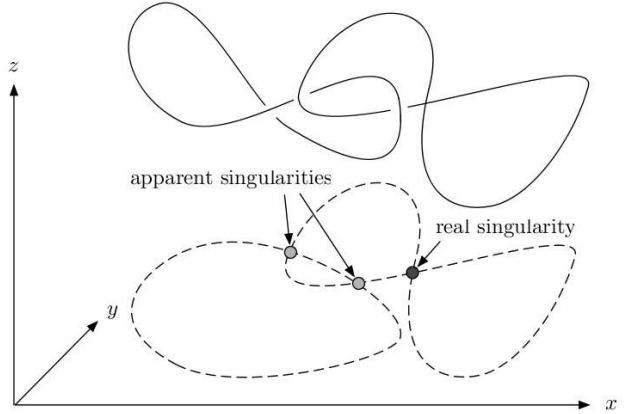}
  \caption{Apparent and real singularities.}
  \label{fig:sweep_topology_singularities}
\end{figure}

\begin{definition}
  {\tmstrong{[Apparent singularity, Real singularity]}}
  
  We call:
  \begin{enumeratenumeric}
    \item Apparent singularities: the singularities of the projected
      curve $\mathcal{D}= \Pi_z (\mathcal{C}_{\mathbbm{R}})$ with at
      least two points as inverse-images (see figure
      \ref{fig:sweep_topology_singularities}).
    
    \item Real singularities: the singularities of the projected curve
      $\mathcal{D}= \Pi_z (\mathcal{C}_{\mathbbm{R}})$ with exactly
      one point as inverse-image (see figure
      \ref{fig:sweep_topology_singularities}).
  \end{enumeratenumeric}
\end{definition}

\begin{definition}
  {\tmstrong{[Generic position]}}
  
  The curve $\mathcal{C}_{\mathbbm{R}}$ is in generic position with respect to
  the $(x, y)$-plane if and only if
  \begin{enumeratenumeric}
    \item $\mathcal{C}_{\mathbbm{R}}$ is in pseudo-generic position with
    respect to the $(x, y)$-plane,
    
    \item $\mathcal{D}= \Pi_z (\mathcal{C}_{\mathbbm{R}})$ is in generic
    position (as a plane algebraic curve)  with respect to
  the $x$-direction,
    
    \item any apparent singularity of $\mathcal{D}= \Pi_z
    (\mathcal{C}_{\mathbbm{R}})$ is a node.
  \end{enumeratenumeric}
\end{definition}
This notion of genericity also appears in a slightly more restrictive form in {\cite{AS}}.

The aim of the next section is to give an algorithm to certify the third point
of the previous definition of generic position. We give also in this section
an effective way to distinguish the real singularities from the apparent ones.

\subsection{Distinguish real singularities and apparent singularities}

In this section, we suppose that $\mathcal{C}_{\mathbbm{R}}$ is in
{\tmstrong{pseudo-generic position}} and $\mathcal{D}= \Pi_z
(\mathcal{C}_{\mathbbm{R}})$ is in {\tmstrong{generic position as a
{\tmstrong{plane}}}} algebraic curve.

Let $(\Gamma_j (X))_{j \in \{ 1,\ldots, n\}}$ be the sequence of $\Gamma$ polynomials associated to the plane curve $\mathcal{D}$ and $(\beta_j (X))_{j \in
\{1,\ldots, n \}}$ be the sequence of associated rational parametrization (see (\ref{betaa}) ).\\
Let $\left( \tmop{Sr}_j (X, Y, Z) \right)_{j \in \{ 0,\ldots, m
\}}$ be the subresultant sequence associated to $P_1, P_2 \in
\mathbbm{Q}[X, Y, Z]$. For any $\left. (k, i) \in \{ 1,\ldots, m
\right\} \times \{ 0,\ldots, k - 1 \left\} \right.$ let,
$ 
R_{k, i} (X, Y)$ be the polynomial\\$k (k - i) \tmop{sr}_{k, i} (X, Y)
   \tmop{sr}^{_{}}_{k, k} (X, Y) - (i + 1) \tmop{sr}_{k, k - 1} (X, Y)
   \tmop{sr}_{k, i + 1} (X, Y).
$
\begin{Lemma}
  Let $(a, b) \in \mathbbm{R}^2$ such that $\tmop{sr}_{k, k} (a, b)\neq 0$, the polynomial $\tmop{Sr}_k (a, b, Z) =
  \sum_{i = 0}^k \tmop{sr}_{k, i} (a, b) Z^i \in \mathbbm{R}[Z]$ has one and
  only one root if and only if $\forall i \in \{ 0,\ldots, k - 1
  \left\} \right.  R_{k, i} (a, b) = 0$.
\end{Lemma}
For any $j \in \{ 1,\ldots, n \left\} \right.$ we define the
sequences\\$(u_k (X))_{k \in \{1,\ldots, j \}}$ and $(v_k (X))_{k \in
\{ 2,\ldots, j \}}$ by\\
$\begin{array}{l}
\ \ u_1 (X) \assign \tmop{gcd} (\Gamma_j (X), \tmop{sr}_{1, 1} (X, \beta_j
(X))),\\
\ \ u_k (X) := \tmop{gcd} (\tmop{sr}_{k, k} (X, \beta_j (X)), u_{k -
1} (X))\\
\ \ v_k (X) \assign \tmop{quo} (u_{k - 1} (X), u_k (X)).
\end{array}$\\
For $k \in \{ 2,\ldots,j \}$ and $i \in \{ 0, k-1 \}$, we define $(w_{k, i} (X))$ by\\
$\begin{array}{l}
\ \ w_{k, 0} (X) \assign v_k (X),\\
\ \ w_{k, i + 1} (X) \assign \tmop{gcd}(R_{k, i} (X, \beta_j (X)), w_{k, i} (X)).
\end{array}
$

\begin{Theorem}
  \label{appsing}For any $j \in \{ 1,\ldots, n \left\} \right.$, let
  $( \Gamma_{j, k} (X) \left) \right. _{\left. k \in \right\{ 1,\ldots, j
  \}}$ and $\left( \chi_{j, k} (X) \left)_{} \right. \right.$ be the
  sequences defined by the following relations
  
  $\Gamma_{j, 1} (X) = \tmop{quo} (\Gamma_j (X), u_1 (X))$ and $\Gamma_{j, k}
  (X) : = w_{k, k} (X)$. $\chi_{j, k} (X) \assign \tmop{quo} (w_{k, 0} (X),
  \Gamma_{j, k} (X))$.
  \begin{enumeratenumeric}
    \item For any root $\alpha$ of $\Gamma_{j, k} (X)$, the x-critical fiber
    $(\alpha, \beta_j (\alpha))$ contain only the point $(\alpha, \beta_j
    (\alpha), \gamma_j (\alpha))$ with \ $\gamma_j (\alpha) \assign -
    \frac{\tmop{sr}_{k, k - 1} (\alpha, \beta_j (\alpha))}{k 
    \tmop{sr}_{k, k} (\alpha, \beta_j (\alpha))}$, so $(\alpha, \beta_j
    (\alpha))$ is a real singularity.
    
    \item For any root $\alpha$ of \ $\chi_{j, k} (X)$, $(\alpha, \beta_j
    (\alpha))$ is an apparent singularity.
    
    \item $\mathcal{C}_{\mathbbm{R}}$ is in generic position if and only if
    for any $(j, k) \in \{ 2,\ldots, n \} \times
    \{ 2,\ldots j \}$ $\chi_{j, k} (X) = 1$.
  \end{enumeratenumeric}
\end{Theorem}

\begin{proof}
  
  \begin{enumeratenumeric}
    \item Let $\alpha$ be a root of \ $\Gamma_{j, k} (X):= w_{k, k} (X) = \gcd (R_{k, k
    - 1} (X, \beta_j (X)), w_{k, k - 1} (X))$. Then $w_{k, k - 1}
    (\alpha) = R_{k, k-1} (\alpha, \beta_j (\alpha)) = 0.$\\
    $w_{k, k - 1} (X) \assign \tmop{gcd} (R_{k,
    k - 2} (X, \beta_j (X)), w_{k, k - 2} (X))$, so\\$ w_{k, k - 2}
    (\alpha) = R_{k, k-2} (\alpha, \beta_j (\alpha)) = 0$.\\
    By induction, using the same argument, it comes that for $i$ from 0 to $(k - 1)$, $w_{k,
    i} (\alpha) =R_{k, i} (\alpha, \beta_j (\alpha)) = 0.$\\
    $w_{k, 0} (X) \assign v_k (X)$, so $v_k (\alpha) = 0$.
    Knowing that $v_k (X) \assign \tmop{quo} (u_{k - 1} (X), u_k (X))$; \
    $u_k (X)$ and $u_{k - 1} (X)$ are square free, then $u_{k - 1} (\alpha) =
    0$ and $u_k (\alpha) \neq 0$. Knowing that $u_k (X) = \tmop{gcd}
    (\tmop{sr}_{k, k} (X, \beta_j (X)), u_{k - 1} (X))$, then\\$
    \tmop{sr}_{k, k} (\alpha, \beta_j (\alpha)) \neq 0$.\\
    $u_{k - 1} (X) = \tmop{gcd} (\tmop{sr}_{k - 1, k - 1} (X, \beta_j (X)),
    u_{k - 2} (X))$ and\\$u_{k - 1} (\alpha) = 0$, so $\tmop{sr}_{k -
    1, k - 1} (\alpha, \beta_j (\alpha)) = u_{k - 2} (\alpha) = 0$.\\
    By induction, using the same argument, it comes that for $i$ from 0 to $k - 1$
    $\tmop{sr}_{i, i} (\alpha, \beta_j (\alpha)) = 0.$\\
    For $i$ from 0 to $k - 1$ $\tmop{sr}_{i, i} (\alpha, \beta_j (\alpha)) =
    0$ and\\$\tmop{sr}_{k, k} (\alpha, \beta_j (\alpha)) \neq 0$, so by the
    fundamental theorem of subresultants,\\
    gcd$(P_1 (\alpha, \beta_j (\alpha), Z), P_2 (\alpha, \beta_j (\alpha),
    Z)) = \tmop{Sr}_k (\alpha, \beta_j (\alpha), Z)\\= \sum_{i = 0}^k
    \tmop{sr}_{k, i} (\alpha, \beta_j (\alpha)) Z^i$. Knowing that\\ gcd$(P_1 (\alpha, \beta_j (\alpha), Z), P_2 (\alpha,
    \beta_j (\alpha), Z)) = \tmop{Sr}_k (\alpha, \beta_j (\alpha), Z)\\=
    \sum_{i = 0}^k \tmop{sr}_{k, i} (\alpha, \beta_j (\alpha)) Z^i$ and \ for
    $i$ from 0 to $(k - 1)$,\\ $R_{k, i} (\alpha, \beta_j (\alpha)) = 0$ then by
    the previous lemma the polynomial gcd$(P_1 (\alpha, \beta_j (\alpha), Z),
    P_2 (\alpha, \beta_j (\alpha), Z)$ have only one root $\gamma_j
    (\alpha) \assign - \frac{\tmop{sr}_{k, k - 1} (\alpha, \beta_j
    (\alpha))}{k \times \tmop{sr}_{k, k} (\alpha, \beta_j (\alpha))}$.
    
    \item Let $\alpha$ be a root of the polynomial \\$ \chi_{j, k} (X)\assign \tmop{quo} (w_{k, 0} (X), \Gamma_{j,
    k} (X))$. Then $w_{k, 0} (\alpha) = 0$ and $\Gamma_{j, k} (\alpha) = w_{k, k} (\alpha)
    \neq 0$ because $w_{k, 0} (X)$ and $\Gamma_{j, k} (X)$ are square free. For $i$ from 0 to $k - 1$, knowing that\\$w_{k, i + 1} (X) \assign
    \tmop{gcd} (R_{k, i} (X, \beta_j (X)), w_{k, i} (X))$, $w_{k, 0} (\alpha)
    = 0$ and $w_{k, k} (\alpha) \neq 0$, then it exist $i \in \{ 0,\ldots, k
    - 1 \left\} \right.$ such that $R_{k, i} (\alpha, \beta_j
    (\alpha)) \neq 0$. So by the previous lemma the polynomial $\tmop{Sr}_k
    (\alpha, \beta_j (\alpha), Z) = \sum_{i = 0}^k \tmop{sr}_{k, i} (\alpha,
    \beta_j (\alpha)) Z^i$ has at least two distinct roots.\\
    By definition $w_{k, 0} (X) \assign v_k (X)$, so $v_k (\alpha) = 0$.
    Knowing that $v_k (X) \assign \tmop{quo} (u_{k - 1} (X), u_k (X))$;\
    $u_k (X)$ and $u_{k - 1} (X)$ are squarefree, then $u_{k - 1} (\alpha) =
    0$ and $u_k (\alpha) \neq 0$.\\
    $u_{k - 1} (\alpha) = 0$, $u_k (\alpha) \neq 0$ and\\$u_k (X) = \tmop{gcd}
    (\tmop{sr}_{k, k} (X, \beta_j (X)), u_{k - 1} (X))$ imply\\$
    \tmop{sr}_{k, k} (\alpha, \beta_j (\alpha)) \neq 0$.\\
    $u_{k - 1} (X) = \tmop{gcd} (\tmop{sr}_{k - 1, k - 1} (X, \beta_j (X)),
    u_{k - 2} (X))$ and\\ $u_{k - 1} (\alpha) = 0$ imply $\tmop{sr}_{k -
    1, k - 1} (\alpha, \beta_j (\alpha)) = u_{k - 2} (\alpha) = 0$.\\
    By induction, using the same argument it comes that for $i$ from 0 to $(k - 1)$
    $\tmop{sr}_{i, i} (\alpha, \beta_j (\alpha)) = 0.$\\
    For $i$ from 0 to $(k - 1)$ $\tmop{sr}_{i, i} (\alpha, \beta_j (\alpha)) =
    0$ and\\$\tmop{sr}_{k, k} (\alpha, \beta_j (\alpha)) \neq 0$, so by the
    fundamental theorem of subresultants\\
    gcd$(P_1 (\alpha, \beta_j (\alpha), Z), P_2 (\alpha, \beta_j (\alpha),
    Z)) = \tmop{Sr}_k (\alpha, \beta_j (\alpha), Z) = \sum_{i = 0}^k
    \tmop{sr}_{k, i} (\alpha, \beta_j (\alpha)) Z^i$.\\
    gcd$(P_1 (\alpha, \beta_j (\alpha), Z), P_2 (\alpha, \beta_j (\alpha),
    Z)) = \tmop{Sr}_k (\alpha, \beta_j (\alpha), Z)$ and $\tmop{Sr}_k (\alpha,
    \beta_j (\alpha), Z)$ has at least two distinct roots imply that $(\alpha,
    \beta_j (\alpha))$ is an apparent singularity.
    
    \item $\mathcal{C}_{\mathbbm{R}}$ is in generic position if and only if
    any apparent singularity of $\mathcal{D}= \Pi_z
    (\mathcal{C}_{\mathbbm{R}})$ is a node. Knowing that the apparent
    singularities of $\mathcal{D}$ which are nodes are exactly those with a
    root of $\chi_{1, 2} (X)$ as $x$-coordinate, so
    $\mathcal{C}_{\mathbbm{R}}$ is in generic position if and only if for any\\
    $\left. (j, k) \in \left\{ 2,\ldots, n \right\} \times \{
    2,\ldots, j \right\}$, $\chi_{j, k} (X) = 1$.
  \end{enumeratenumeric}
\end{proof}

\subsection{Lifting and connection phase}

In this section, we suppose that $\mathcal{C}_{\mathbbm{R}}$ is in
{\tmstrong{generic position}} that means that $\mathcal{C}_{\mathbbm{R}}$ is
in {\tmstrong{pseudo-generic position}}, $\mathcal{D}= \Pi_z
(\mathcal{C}_{\mathbbm{R}})$ is in {\tmstrong{generic position as a
{\tmstrong{plane}}}} algebraic curve and any {\tmstrong{apparent singularity}}
of $\mathcal{D}= \Pi_z (\mathcal{C}_{\mathbbm{R}})$ is a {\tmstrong{node}}.\\
To compute the topology of $\mathcal{C}_{\mathbbm{R}}$ we first compute the
topology of its projection on the $(x, y)$-plane and in second we lift the computed topology.

As mentioned in section 2, to compute the topology
of a plane algebraic curve in generic position, we need to compute its critical
fibers and one regular fiber between two critical ones. So to obtain the
topology of $\mathcal{C}_{\mathbbm{R}}$ we just need to lift the critical and
regular fibers of $\mathcal{D}= \Pi_z (\mathcal{C}_{\mathbbm{R}})$.

Here after we explain how this lifting can be done without any supplementary
computation for the regular fibers and the real critical fibers. And for the
special case of the apparent singular fibers, we present a new approach for the lifting and the connections.

\subsubsection{Lifting of the regular points of $\mathcal{D}= \Pi_z
(\mathcal{C}_{\mathbbm{R}})$}

The lifting of the regular fibers of $\mathcal{D}= \Pi_z
(\mathcal{C}_{\mathbbm{R}})$ is done by using the rational parametrizations given in Proposition \ref{proplifting}.
%

\subsubsection{Lifting of the real singularities of $\mathcal{D}= \Pi_z
(\mathcal{C}_{\mathbbm{R}})$}

The lifting of the real singularities of $\mathcal{D}= \Pi_z
(\mathcal{C}_{\mathbbm{R}})$ is done by using the rational parametrizations
given by 1. of Theorem \ref{appsing}.

\subsubsection{Connection between real singularities and regular points}
\begin{figure}
  \centering
  \includegraphics[width=\linewidth]{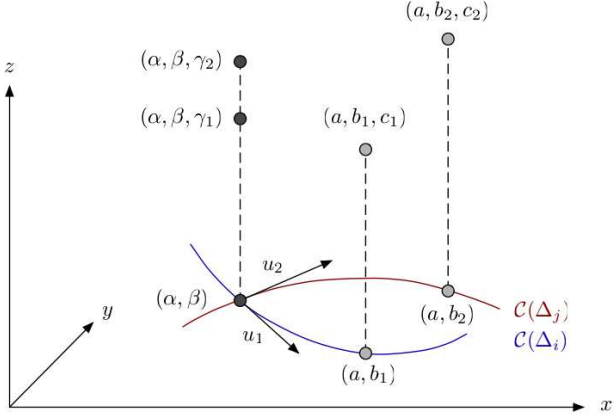}
  \caption{Connection between real singularities and regular points.}
  \label{fig:lifting_projected_curves}
\end{figure}

For a space curve in pseudo-generic position, the connections between real singularities and regular points are exactly those obtained on the projected curve using Grandine's sweeping algorithm {\cite{GVN}} (see figure \ref{fig:lifting_projected_curves}).

\subsubsection{Lifting of the apparent singularities}

\begin{figure}
  \centering
  \includegraphics[width=\linewidth]{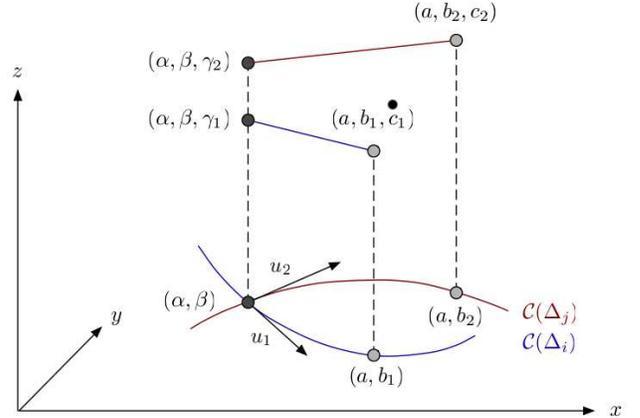}
  \caption{Lifting of an apparent singularity.}
  \label{fig:lifting_apparent_singularity_pre}
\end{figure}

The lifting of the topology around an apparent singularity is a little more
complex. Above an apparent singularity of $\mathcal{D}= \Pi_z
(\mathcal{C}_{\mathbbm{R}})$ we have firstly to compute the $z$-coordinates
and secondly to decide which of the two branches pass over the other (see figure \ref{fig:lifting_apparent_singularity_pre}). We solve these problems by analyzing the situation at an apparent singularity.

According to Theorem \ref{thmpseudogenpos} (2.), $\mathcal{D}= \Pi_z
(\mathcal{C}_{\mathbbm{R}}) = \underset{i = 1}{\overset{m}{\bigcup}}
\mathcal{C}(\Delta_i),$ so an apparent singularity is a cross point of a branch of  $\mathcal{C}(\Delta_i)$ and a branch of 
$\mathcal{C}(  \Delta_j)$ with $i, j \in \{ 1,\ldots, m \}$. So we have the
following proposition.

\begin{Proposition}\label{toto}
  If $(\alpha, \beta)$ is an apparent singularity of $\mathcal{D}$ such that
  $\Delta_i (\alpha, \beta) = \Delta_j (\alpha, \beta) = 0$, then the degree
  of the polynomial $\gcd (P_1 (\alpha, \beta, Z), P_2 (\alpha, \beta, Z)) \in
  \mathbbm{R}[Z]$ will be $(i + j)$. 
\end{Proposition}

Let $(\alpha, \beta)$ be an apparent singularity of $\mathcal{D}$ such that\\
$\Delta_i (\alpha, \beta) = \Delta_j (\alpha, \beta) = 0$ and $\gamma_1$,
$\gamma_2$ the corresponding $z$-coordinates. So by Proposition \ref{toto} and Proposition
\ref{propsubres} \\$\tmop{sr}_{0, 0} (\alpha, \beta) = \ldots = \tmop{sr}_{i,
i} (\alpha, \beta) = \ldots = \tmop{sr}_{j, j} (\alpha, \beta) = . \ldots =
\tmop{sr}_{i + j - 1, i + j - 1} (\alpha, \beta) = 0$.\\
By Proposition \ref{proplifting}, for any $(a, b, c) \in
\mathcal{C}_{\mathbbm{R}}$ such that $\Delta_i (a, b) = 0$ and $\tmop{sr}_{i,
i} (a, b) \neq 0$ we have $c = - \frac{\tmop{sr}_{i, i - 1} (a, b)}{i
\tmop{sr}_{i, i} (a, b)}$. So the function $(x, y) \longmapsto Z_i \assign -
\frac{\tmop{sr}_{i, i - 1} (x, y)}{i \tmop{sr}_{i, i} (x, y)}$ gives the
$z$-coordinate of any $(a, b, c) \in \mathcal{C}_{\mathbbm{R}}$ such that
$\Delta_i (a, b) = 0$ and $\tmop{sr}_{i, i} (a, b) \neq 0$.\\
$\Delta_i (\alpha, \beta) = 0$ but $\tmop{sr}_{i, i} (\alpha, \beta) = 0$, so
the function $Z_i$ is not defined on $(\alpha, \beta)$. The solution comes
from the fact that the function $Z_i$ is continuously extensible on $(\alpha, \beta)$. Let
$u_1$ be the slope of the tangent line of $\mathcal{C}(\Delta_i)$ at
$(\alpha, \beta)$ and $t \in \mathbbm{R}^{\ast}$. Let $\gamma_i (t) \assign
Z_i (\alpha, \beta + tu_1) = - \frac{\tmop{sr}_{i, i - 1} (\alpha, \beta +
tu_1)}{i \tmop{sr}_{i, i} (\alpha, \beta + tu_1)}$. Knowing that the algebraic
curve $\mathcal{C}_{\mathbbm{R}}$ hasn't any discontinuity, it comes $\lim_{t
\rightarrow 0^+} \gamma_i (t) = \lim_{t \rightarrow 0^-} \gamma_i (t) =
\gamma_1$. By the same arguments, if we denote $u_2$ the slope of the
tangent line of \ $\mathcal{C}(\Delta_j)$ at $(\alpha, \beta)$ and $\gamma_j
(t) \assign Z_j (\alpha, \beta + tu_2) = - \frac{\tmop{sr}_{j, j - 1} (\alpha,
\beta + tu_2)}{j \tmop{sr}_{j, j} (\alpha, \beta + tu_2)}$, then $\lim_{t
\rightarrow 0^+} \gamma_j (t) = \lim_{t \rightarrow 0^-} \gamma_j (t) =
\gamma_2$. The values $u_1, u_2, \gamma_1$ and $\gamma_2$ are computed using Taylor formulas and certified numerical approximations. 

\begin{figure}
  \centering
  \includegraphics[width=\linewidth]{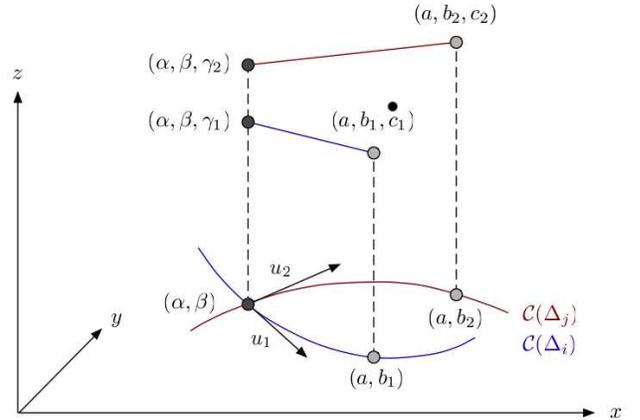}
  \caption{Connection above an apparent singularity.}
  \label{fig:lifting_apparent_singularity}
\end{figure}

Now it remains to decide which of the two branches pass over the other. This
problem is equivalent to the problem of deciding the connection around an
apparent singularity. Let $(a, b_1, c_1)$ and $(a, b_2, c_2)$ the regular
points that we have to connect to $(\alpha, \beta, \gamma_1)$ and $(\alpha,
\beta, \gamma_2)$. The question is which of the points $(a, b_1, c_1)$ and $(a,
b_2, c_2)$ will be connected to $(\alpha, \beta, \gamma_1)$ and the other to
$(\alpha, \beta, \gamma_2)$ (see figure \ref{fig:lifting_apparent_singularity_pre})? In {\cite{AS}} Alc\'azar and Sendra give a solution using a second projection of the space curve but it costs a computation of a Sturm Habicht sequence of $P_1$ and $P_2$. Our solution does not use any supplementary computation. It comes from the fact that $\gamma_1$
is associated to $u_1$ and $\gamma_2$ to $u_2$. Knowing that $u_1$ is the
slope of the tangent line of $\mathcal{C}(\Delta_i)$ at $(\alpha, \beta)$
and $u_2$ the slope of the tangent line of \ $\mathcal{C}(\Delta_j)$ at
$(\alpha, \beta)$, so $(\alpha, \beta, \gamma_1)$ will be connected to $(a,
b_1, c_1)$ if $(a, b_1)$ is on the branch associated to $u_1$. If $(a, b_1)$
is not on the branch associated to $u_1$, then $(a, b_1)$ is on the branch associated to $u_2$, so $(\alpha, \beta, \gamma_2)$
will be connected to $(a, b_1, c_1)$ (see figure \ref{fig:lifting_apparent_singularity}).

\begin{figure*}[!htb]
\centering
\begin{tabular*}{\linewidth}{@{\extracolsep{\fill}}|c|c|c|c|}
\hline 
Curve & $P_1(x, y, z)$ &$P_2(x, y, z)$ & Time (s)\\
\hline 
1 &   $ x^2+y^2+z^2-1$&$x^2-y^2-z+1 $ & $0.032$ \\
\hline 
2 & $ x^2+y^2+z^2-1$&$x^3+3x^2z+3xz^2+z^3+y^3-xyz-yz^2 $ & $0.659$ \\
\hline
3 & $ (x-2y+2z)^2+y^2+z-1$&$z^3-z-(x-2y+2z)^3+3(x-2y+2z)y^2$ & $2.125$ \\
\hline
4 & $ (x-2y+2z)^2+y^2+z^2-1$&$y^3-(x-2y+2z)^3-(x-2y+2z)yz$ & $1.031$ \\
\hline
5 & $ (x-y+z)^2+y^2+z^2-1$&$y^2-(x-y+z)^2-(x-y+z)z)^2-z^2((x-y+z)^2+y^2)$ & $1.6963$ \\
\hline
6 & $ (x-y+z)^2+y^2+z^2-1$&$((x-y+z)^2+y^2+z^2)^2-4((x-y+z)^2+y^2)$ & $2.228$ \\
\hline
7 & $(x-y+z)^2+y^2-2(x-y+z)$&$((x-y+z)^2+y^2+z^2)^2-4((x-y+z)^2+y^2)$ & 2.875\\
\hline
\end{tabular*}
\caption{Running time of experimentations.}
\label{fig:running_time}
\end{figure*}

\begin{Remark} For a curve in generic position any apparent singularity is a node, so the slopes at an apparent singularity are always distinct that is to say $u_1 \neq u_2$.
\end{Remark}

\section{Implementation, experiments}
A preliminary implementation of our method has been written using the Computer Algebra System Mathemagix. Results are visualized using the Axel\footnote{\tt http://axel.inria.fr} algebraic geometric modeler which allows the manipulation of geometric objects with algebraic representation such as implicit or parametric curves or surfaces.

Since existing methods have no publicly available implementations, table \ref{fig:running_time} only reports our experiments, performed on an Intel(R) Core machine clocked at 2GHz with 1GB RAM.

\begin{figure}[t!]
  \centering
  \includegraphics[width=\linewidth]{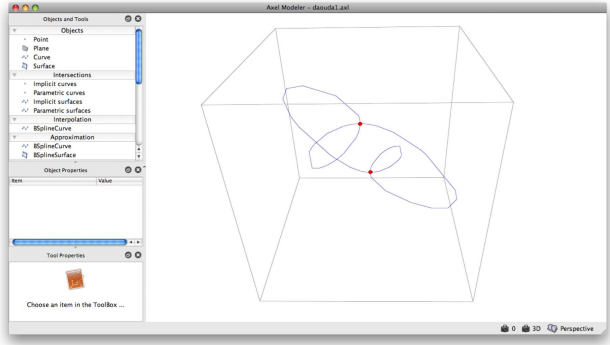}
  \caption{Computed topology of curve 2 of table \ref{fig:running_time}.} 
  \label{fig:SpaceCurve2}
\end{figure}

\begin{figure}[t!]
  \centering
  \includegraphics[width=\linewidth]{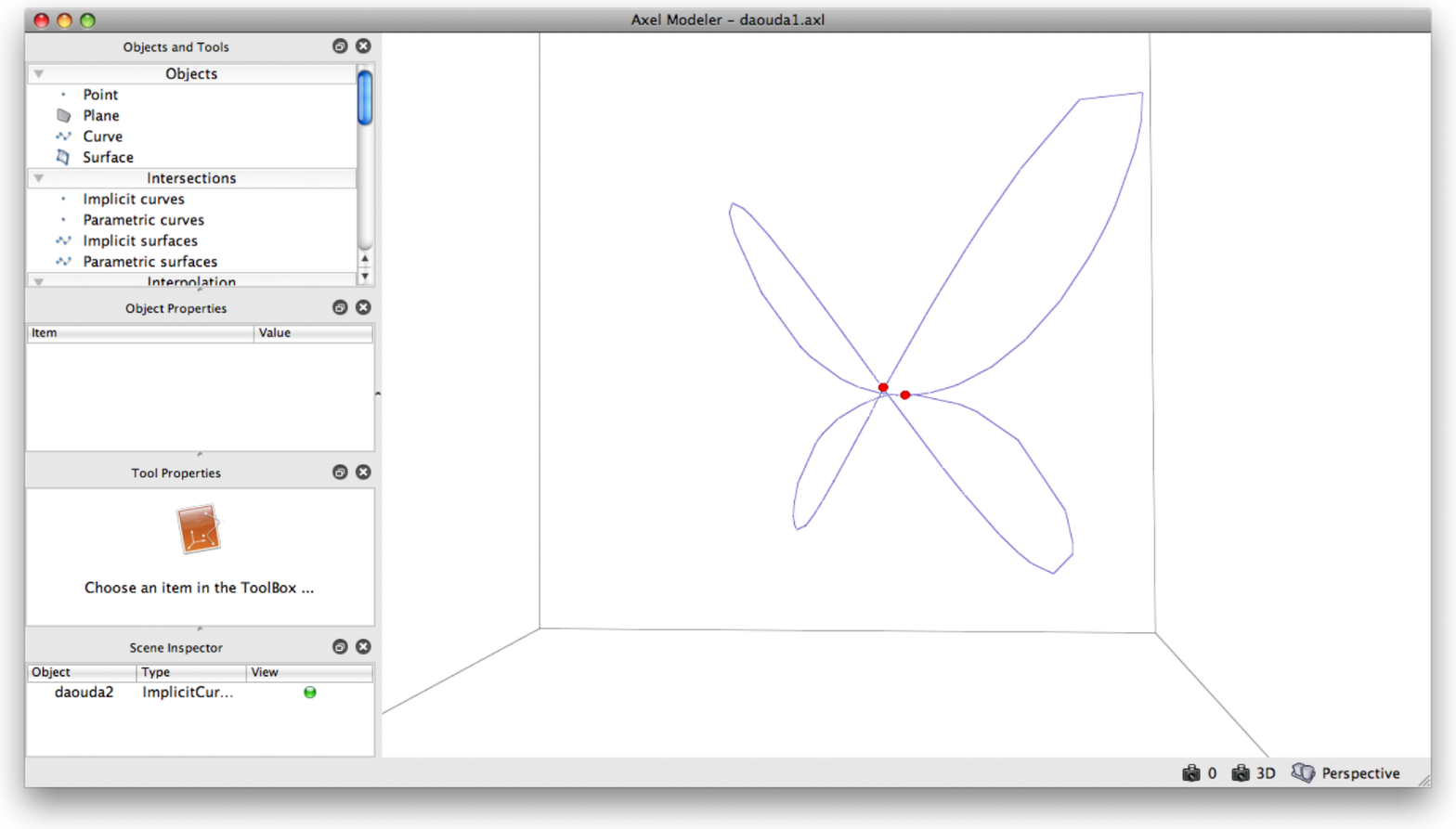}
  \caption{Computed topology of curve 7 of table \ref{fig:running_time}.}
  \label{fig:SpaceCurve1}
\end{figure}

\end{document}